\documentclass[a4paper,11pt]{article}%
\usepackage{graphicx}
\usepackage{geometry}
\usepackage{amsmath}
\usepackage{amssymb}
\usepackage{amsfonts}
\usepackage{amsthm,amscd}%
\providecommand{\U}[1]{\protect\rule{.1in}{.1in}}
\def\figurename{Figure}
\makeatletter
\renewcommand{\fnum@figure}[1]{\figurename~\thefigure.}
\makeatother
\def\tablename{Table}
\makeatletter
\renewcommand{\fnum@table}[1]{\tablename~\thetable.}
\makeatother
\def \bop {\noindent\textbf{Proof. }}
\def \eop {\hbox{}\nobreak\hfill
\vrule width 2mm height 2mm depth 0mm
\par \goodbreak \smallskip}
\newtheorem{theorem}{Theorem}[section]
\newtheorem{lemma}[theorem]{Lemma}

\newtheorem{proposition}[theorem]{Proposition}
\theoremstyle{definition}
\newtheorem{definition}[theorem]{Definition}

\theoremstyle{remark}
\newtheorem{remark}[theorem]{Remark}
\numberwithin{equation}{section}

\begin{document}

\title{Stability of Mc Kean-Vlasov stochastic differential \\equations and applications}
\author{{Khaled Bahlali}\thanks{ Laboratoire IMATH, Universit\'{e} du Sud-Toulon-Var,
B.P 20132, 83957 La Garde Cedex 05, France. \textit{(E-mail:
bahlali@univ-tln.fr)}}
\and {Mohamed Amine Mezerdi }\thanks{ Laboratory of Applied Mathematics, University
of Biskra, Po. Box 145, Biskra (07000), Algeria. \textit{(E-mail:
mohamed@live.com)}}
\and Brahim Mezerdi\thanks{ King Fahd University of Petroleum and Minerals,
Department of Mathematics and Statistics, P.O. Box 1916, Dhahran 31261, Saudi
Arabia. (\textit{E-mail: brahim.mezerdi@kfupm.edu.sa)}}}
\maketitle
\date{}

\begin{abstract}
We consider Mc Kean-Vlasov stochastic differential equations (MVSDEs), which
are SDEs where the drift and diffusion coefficients depend not only on the
state of the unknown process but also on its probability distribution. This
type of SDEs was studied in statistical physics and represents the natural
setting for stochastic mean-field games. We will first discuss questions of
existence and uniqueness of solutions under an Osgood type condition improving
the well known Lipschitz case. Then we derive various stability properties
with respect to initial data, coefficients and driving processes, generalizing
known results for classical SDEs. Finally, we establish a result on the
approximation of \ the solution of a MVSDE associated to a relaxed control by
the solutions of the same equation associated to strict controls. As a
consequence, we show that the relaxed and strict control problems have the
same value function. This last property improves known results proved for a
special class of MVSDEs, where the dependence on the distribution was made via
a linear functional.



\end{abstract}

\section{Introduction}

We will investigate some properties of a particular class of stochastic
differential equations (SDE), called Mc Kean-Vlasov stochastic differential
equations (MVSDE) or mean-field stochastic differential equations. These are
SDEs described by

\begin{center}
$\left\{
\begin{array}
[c]{l}%
dX_{t}=b(t,X_{t},\mathbb{P}_{X_{t}})ds+\sigma(t,X_{t},\mathbb{P}_{X_{t}%
})dB_{s}\\
X_{0}=x,
\end{array}
\right.  $
\end{center}

where $b$ is the drift, $\sigma$ is the diffusion coefficient and $\left(
B_{t}\right)  $ is a Brownian motion. For this type of equations the drift and
diffusion coefficient depend not only on the state variable $X_{t},$ but also
on its marginal distribution $\mathbb{P}_{X_{t}}$. This fact brings a non
trivial additional difficulty compared to classical It\^{o} SDEs. The
solutions of such equation are known in the literature as non linear diffusions.

MVSDEs were first studied in statistical physics by M. Kac \cite{Ka}, as a
stochastic counterpart for the Vlasov equation of plasma \cite{Vla}. The
probabilistic study of such equation has been performed by H.P. Mc Kean
\cite{McK}, see \cite{Sn} for an introduction to this research field. These
equations were obtained as limits of some weakly interacting particle systems
as the number of particles tends to infinity. This convergence property is
called in the literature as the propagation of chaos. The MVSDE, represents in
some sense the average behavior of the infinite number of particles. One can
refer to \cite{CarDel, Gra, JMW} for details on the existence and uniqueness
of solutions for such SDEs, see also \cite{BuDjLiPe, BuLiPe} for the case of
Mc Kean Vlasov backward stochastic differential equations (MVBSDE). Existence
and uniqueness with less regularity on the coefficients have been established
in \cite{Ch, ChFr, Chia, HSS, MiVe, Sheu}. Recently there has been a renewed
interest for MVSDEs, in the context of mean-field games (MFG) theory,
introduced independently by P.L. Lions and J.M. Lasry \cite{LasLio} and Huang,
Malham\'{e} Caines \cite{HMC} in 2006. MFG theory has been introduced to solve
the problem of existence of an approximate Nash equilibrium for differential
games, with a large number of players (see \cite{Ben}). Since the earlier
papers, MFG theory and mean-field control theory has raised a lot of interest,
motivated by applications to various fields such as game theory, mathematical
finance, communications networks and management of oil ressources. One can
refer to the most recent and updated reference on the subject \cite{CarDel}
and the complete bibliographical list therein.

Our main objective in this paper is to study some properties of such equations
such as existence, uniqueness, and stability properties. In particular, we
prove an existence and uniqueness theorem for a class of MVSDEs under Osgood
type condition on the coefficients, improving the well known globally
Lipschitz case.\ It is well known that stability properties of deterministic
or stochastic dynamical systems are crucial in the study of such systems. It
means that the trajectories do not change too much under small perturbations.
We study stability with respect to initial conditions, coefficients and
driving processes, which are continuous martingales and bounded variation
processes. These properties will be investigated under Lipschitz condition
with respect to the state variable and the distribution and generalize known
properties for classical It\^{o} SDEs, see \cite{BMO, IW}. Furthermore, we
prove that in the context of stochastic control of systems driven by MVSDEs,
the relaxed and strict control problems have the same value function. As it is
well known when the Filipov type convexity condition is not fulfilled, there
is no mean to prove the existence of a strict control. The idea is then to
embedd the usual strict controls into the set of measure valued controls,
called relaxed controls, which enjoys good compactness properties. So for the
relaxed control to be a true extension of the initial problem, the value
functions of both control problems must be the same. Under the Lipschitz
condition we prove that the value functions are equal. Note that this result
extends to general Mc Kean Vlasov equations known results \cite{BMM1, BMM2}
established for a special class of MVSDEs, where the dependence of the
coefficient on the distribution variable is made via a linear form of the distribution.

\section{Formulation of the problem and preliminary results}

\subsection{Assumptions}

Let $(\Omega,\mathcal{F},P)$ be a probability space$,$ equipped with a
filtration $\left(  \mathcal{F}_{t}\right)  ,$ satisfying the usual conditions
and $\left(  B_{t}\right)  $ a $d$-dimensional $\left(  \mathcal{F}%
_{t},P\right)  -$Brownian motion. Let us consider the following Mc Kean-Vlasov
stochastic differential equation called also mean-field stochastic
differential equation (MVSDE)

\begin{center}%
\begin{equation}
\left\{
\begin{array}
[c]{l}%
dX_{t}=b(t,X_{t},\mathbb{P}_{X_{t}})ds+\sigma(t,X_{t},\mathbb{P}_{X_{t}%
})dB_{s}\\
X_{0}=x
\end{array}
\right.  \label{MVSDE}%
\end{equation}

\end{center}

Note that for this kind of SDEs, the drift $b$ and diffusion coefficient
$\sigma$ depend not only on the position, but also on the marginal
distribution of the solution.

The following assumption will be considered throughout this paper.

Let us denote $\mathcal{P}_{2}(\mathbb{R}^{d})$ the space of probability
measures with finite second order moment. That is for each $\mu\in
\mathcal{P}_{2}(\mathbb{R}^{d})$ $%
{\displaystyle\int}
\left\vert x\right\vert ^{2}\mu(dx)<+\infty.$

(\textbf{H}$_{\mathbf{1}}$\textbf{)} Assume that

\begin{center}
$%
\begin{array}
[c]{c}%
b:[0,T]\times\mathbb{R}^{d}\times\mathcal{P}_{2}(\mathbb{R}^{d}%
)\longrightarrow\mathbb{R}^{d}\\
\sigma:[0,T]\times\mathbb{R}^{d}\times\mathcal{P}_{2}(\mathbb{R}%
^{d})\longrightarrow\mathbb{R}^{d}\otimes\mathbb{R}^{d}%
\end{array}
$
\end{center}

are Borel measurable functions and there exist $C>0$ such that for every
$(t,x,\mu)\in\lbrack0,T]\times\mathbb{R}^{d}\times\mathcal{P}_{2}%
(\mathbb{R}^{d}):$

\begin{center}
$|b(t,x,\mu)|+|\sigma(t,x,\mu)|\leq C\left(  1+\left\vert x\right\vert
\right)  $
\end{center}

(\textbf{H}$_{\mathbf{2}}$\textbf{)} There exist $L>0$ such that for any
$t\in\lbrack0,T],x,$ $x^{\prime}\in\mathbb{R}^{d}$ and $\mu,$ $\mu^{\prime}%
\in\mathcal{P}_{2}(\mathbb{R}^{d}),$

\begin{center}
$|b(t,x,\mu)-b(t,x^{\prime},\mu^{\prime})|\leq L[|x-x^{\prime}|+W_{2}(\mu
,\mu^{\prime})]$

$|\sigma(t,x,\mu)-\sigma(t,x^{\prime},\mu^{\prime})|\leq L[|x-x^{\prime
}|+W_{2}(\mu,\mu^{\prime})]$
\end{center}

where $W_{2}$ denotes the 2-Wasserstein metric.

\subsection{Wasserstein metric}

Let $\mathcal{P}(\mathbb{R}^{d})$ be the space of probability measures on
$\mathbb{R}^{d}$ and for any $p>1$, denote by $\mathcal{P}_{p}(\mathbb{R}%
^{d})$ the subspace of $\mathcal{P}(\mathbb{R}^{d})$ of the probability
measures with finite moment of order $p$.

For $\mu,\nu\in\mathcal{P}_{p}(\mathbb{R}^{d}),$ define the $p$-Wasserstein
distance $W_{p}(\mu,\nu)$ by:
\[
W_{p}(\mu,\nu)=\inf_{\pi\in\Pi(\mu,\nu)}[\int_{E\times E}\left\vert
x-y\right\vert ^{p}d\pi(x,y)]^{1/p}%
\]

where $\Pi(\mu,\nu)$ denotes the set of probability measures on $\mathbb{R}%
^{d}\times\mathbb{R}^{d}$ whose first and second marginals are respectively
$\mu$ and $\nu$.

In the case $\mu=\mathbb{P}_{X}$ and $\nu=\mathbb{P}_{Y}$ are the laws of
$\mathbb{R}^{d}$-valued random variable $X$ and $Y$ of order $p$, then
\[
W_{p}(\mu,\nu)^{p}\leq\mathbb{E[}\left\vert X-Y\right\vert ^{p}].
\]

Indeed%

\begin{align*}
W_{p}(\mu,\nu)  &  =\inf_{\pi\in\Pi(\mu,\nu)}[\int_{E\times E}\left\vert
x-y\right\vert ^{p}d\pi(x,y)]^{1/p}\\
W_{p}(\mu,\nu)^{p}  &  =\inf_{\pi\in\Pi(\mu,\nu)}[\int_{E\times E}\left\vert
x-y\right\vert ^{p}d\pi(x,y)]\\
&  \leq\int_{E\times E}\left\vert x-y\right\vert ^{p}d(\mathbb{P}_{\left(
X,Y\right)  }(x,y)\\
&  =\mathbb{E[}\left\vert X-Y\right\vert ^{p}]
\end{align*}

In the literature the Wasserstein metric is restricted to $W_{2}$ while
$W_{1}$ is often called the Kantorovich-Rubinstein distance because of the
role it plays in optimal transport.

\section{Existence and uniqueness of solutions}

\subsection{The globally Lipschitz case}

The following theorem states that under global Lipschitz condition,
(\ref{MVSDE}) admits a unique solution. Its complete proof is given in
\cite{Sn} for a drift depending linearly on the law of $X_{t}$ that is
$b(t,x,\mu)=%
{\displaystyle\int\limits_{\mathbb{R}^{d}}}
b^{\prime}(t,x,y)\mu(dy)$ and a constant diffusion. The general case as in
(\ref{MVSDE}) is treated in \cite{CarDel} Theorem 4.21 or \cite{JMW}
Proposition 1.2 and is based on a fixed point theorem on the space of
continuous functions with values in $\mathcal{P}_{2}(\mathbb{R}^{d})$. Note
that in \cite{Gra, JMW} the authors consider MVSDEs driven by general L\'{e}vy
process instead of a Brownian motion.

\begin{theorem}
Under assumptions $\mathbf{(H_{1})}$, $\mathbf{(H_{2})}$, (\ref{MVSDE}) admits
a unique solution such that

$E[\sup_{t\leq T}|X_{t}|^{2}]<+\infty$
\end{theorem}

\bop

Let us give the outline of the proof. Let $\mu\in\mathcal{P}_{p}%
(\mathbb{R}^{d})$ be fixed, the classical It\^{o}'s theorem gives the
existence and uniqueness of a solution denote by $\left(  X_{t}^{\mu}\right)
$ satisfying $E[\sup_{t\leq T}|X_{t}^{\mu}|^{2}]<+\infty.$ Now let us consider
the mapping

\begin{center}
$\Psi:\mathcal{C}(\left[  0,T\right]  ,\mathcal{P}_{2}(\mathbb{R}%
^{d}))\longrightarrow\mathcal{C}(\left[  0,T\right]  ,\mathcal{P}%
_{2}(\mathbb{R}^{d}))$

$\mu\longrightarrow\Psi(\mu)=\left(  \mathcal{L}(X_{t}^{\mu})\right)
_{t\geq0},$ the distribution of $X_{t}^{\mu}.$
\end{center}

$\Psi$ is well defined as $X_{t}^{\mu}$ has continuous paths and
$E[\sup_{t\leq T}|X_{t}^{\mu}|^{2}]<+\infty.$

To prove the existence and uniqueness of(\ref{MVSDE}), it is sufficient to
prove that the mapping $\Psi$ has a unique fixed point. By using usual
arguments from stochastic calculus and relation and the property of
Wasserstein metric it is easy to show that:

\begin{center}
$\sup_{t\leq T}W_{2}(\left(  \Psi^{k}(\mu)\right)  _{t},\left(  \Psi^{k}%
(\nu)\right)  _{t})^{2}\leq C\dfrac{T^{k}}{k!}\sup_{t\leq T}W_{2}(\mu_{t}%
,\nu_{t})^{2}$
\end{center}

For large $k,$ $\Psi^{k}$ is a strict contraction which implies that $\Psi$
admits a unique fixed point in the complete metric space $\mathcal{C}(\left[
0,T\right]  ,\mathcal{P}_{2}(\mathbb{R}^{d})).$\eop

The following version MVSDEs is also considered in the control literature%

\begin{equation}
\left\{
\begin{array}
[c]{l}%
dX_{t}=b(t,X_{t},%
{\displaystyle\int}
\varphi(y)\mathbb{P}_{X_{t}}(dy))dt+\sigma(t,X_{t},%
{\displaystyle\int}
\psi(y)\mathbb{P}_{X_{t}}(dy))dW_{t}\\
X_{0}=x
\end{array}
\right.  \label{MVSDE1}%
\end{equation}

where

$\mathbf{(H}_{\mathbf{3}}\mathbf{)}$ $b,$ $\sigma,\varphi$ and $\psi$ are
Borel measurable bounded functions such that $b(t,.,.),$ $\sigma(t,.,.),$
$\varphi$ and $\psi$ are globally lipshitz functions in $\mathbb{R}^{d}%
\times\mathbb{R}^{d}$.

\begin{proposition}
Under assumptions $\mathbf{(H}_{\mathbf{1}}\mathbf{)}$ and $\left(
\mathbf{H}_{\mathbf{3}}\right)  $the MVSDE (\ref{MVSDE1}) has a unique strong
solution. Moreover for each $p>0$ we have $E(\left\vert X_{t}\right\vert
^{p})<+\infty.$
\end{proposition}

\bop

Let us define $\overline{b}(t,x,\mu)$ and $\overline{\sigma}(t,x,\mu)$ on
$\left[  0,T\right]  \times\mathbb{R}^{d}\times\mathcal{P}_{2}(\mathbb{R}%
^{d})$ by

\begin{center}
$\overline{b}(t,x,\mu)=b(.,.,%
{\displaystyle\int}
\varphi(x)d\mu(x),.)$, $\overline{\sigma}(t,x,\mu)=\sigma(t,x,%
{\displaystyle\int}
\psi(x)d\mu(x)).$
\end{center}

According to the last Theorem it is sufficient to check that $\overline{b}$
and $\overline{\sigma}$ are Lipschitz in $\left(  x,\mu\right)  $. Indeed
since the coefficients $b$ and $\sigma$ are Lipschitz continuous in $x,$ then
$\overline{b}$ and $\overline{\sigma}$ are also Lipschitz in $x.$ Moreover one
can verify easily that $\overline{b}$ and $\overline{\sigma}$ are also
Lipshitz continuous in $\mu,$ with respect to the Wasserstein metric

\begin{center}
$%
\begin{array}
[c]{l}%
W_{2}(\mu,\nu)=\inf\left\{  \left(  E^{Q}\left\vert X-Y\right\vert
^{2}\right)  ^{1/2};Q\in\mathcal{P}_{2}(\mathbb{R}^{d}\times\mathbb{R}%
^{d}),\text{ with marginals }\mu,\nu\right\} \\
=\sup\left\{
{\displaystyle\int}
hd\left(  \mu-\nu\right)  ;\text{ }\left\vert h(x)-h(y)\right\vert
\leq\left\vert x-y\right\vert \right\}  ,
\end{array}
$
\end{center}

Note that the second equality is given by the Kantorovich-Rubinstein theorem
\cite{CarDel}. Since the mappings $b$ and $\varphi$ in the the MFSDE are
Lipschitz continuous in $y$ we have

\begin{center}
$%
\begin{array}
[c]{l}%
\left\vert b(.,.,%
{\displaystyle\int}
\varphi(y)d\mu(y),.)-b(.,.,%
{\displaystyle\int}
\varphi(y)d\nu(y),.)\right\vert \\
\leq K\left\vert
{\displaystyle\int}
\varphi(y)d(\mu(y)-\nu(y))\right\vert \\
\leq K^{\prime}.W_{2}\left(  \mu,\nu\right)
\end{array}
$
\end{center}

Therefore $\overline{b}(t,.,.)$ is Lipschitz continuous in the variable
$(x,\mu)\in\mathbb{R}^{d}\times\mathcal{P}_{2}(\mathbb{R}^{d})$ uniformly in
$t\in\left[  0,\left[  T\right]  \right]  $

Similar arguments can be used for $\sigma.$ \eop

\subsection{The uniqueness under Osgood type condition}

In this section we relax the global Lipschitz condition in the state variable.
We will prove the existence and uniqueness of a solution when the coefficients
are globally Lipschitz in the distribution variable and satisfy an Osgood type
condition in the state variable. To be more precise let us consider the
following MVSDE

\begin{center}%
\begin{equation}
\left\{
\begin{array}
[c]{c}%
dX_{t}=b(t,X_{t},\mathbb{P}_{X_{t}})ds+\sigma(t,X_{t})dB_{s}\\
X_{0}=x
\end{array}
\right.  \label{MVSDE2}%
\end{equation}

\end{center}

Assume that $b$ and $\sigma$ are real valued bounded Borel measurable
functions satisfying:

$\mathbf{(H}_{\mathbf{4}}\mathbf{)}$ There exist $C>0,$ such that for every
$x\in\mathbb{R}$ and $\left(  \mu,\nu\right)  \in\mathcal{P}_{1}%
(\mathbb{R})\times\mathcal{P}_{1}(\mathbb{R}):$

\begin{center}
$|b(t,x,\mu)-b(t,x,\nu)|\leq CW_{1}(\mu,\nu)$
\end{center}

$\mathbf{(H}_{\mathbf{5}}\mathbf{)}$There exists a strictly increasing
function $\rho(u)$ on $[0,+\infty)$ such that $\rho(0)=0$ and $\rho^{2}$ is
convex satisfying $%
{\displaystyle\int\limits_{0^{+}}}
\rho^{-2}(u)du=+\infty,$ such that for every $(x,y)\in\mathbb{R}%
\times\mathbb{R}$ and $\mu\in\mathcal{P}_{2}(\mathbb{R}),$ $|\sigma
(t,x)-\sigma(t,y)|\leq\rho(|x-y|).$

$\mathbf{(H}_{\mathbf{6}}\mathbf{)}$ There exists a strictly increasing
function $\kappa(u)$ on $[0,+\infty)$ such that $\kappa(0)=0$ and $\kappa$ is
concave satisfying $%
{\displaystyle\int\limits_{0^{+}}}
\kappa^{-1}(u)du=+\infty,$ such that for every $(x,y)\in\mathbb{R}^{d}%
\times\mathbb{R}^{d}$ and $\mu\in\mathcal{P}_{2}(\mathbb{R}),$ $|b(t,x,\mu
)-b(t,y,\mu)|\leq\kappa(|x-y|).$

In the next Theorem we derive the pathwise uniqueness for (\ref{MVSDE2}) under
an Osgood type condition in the state variable. This result improves
\cite{IW},Theorem 3.2, established for classical It\^{o}'s SDEs and
\cite{CarDel} Theorem 4.21, at least for MVSDEs with a diffusion coefficient
not depending on the distribution variable.

\begin{theorem}
Under assumptions $\mathbf{(H_{4})}-\mathbf{(H_{6})}$, the MVSDE
(\ \ref{MVSDE2}) enjoys the property of pathwise uniqueness.
\end{theorem}

\bop

The following proof is inspired from \cite{CarDel} Theorem 4.21.

Since $%
{\displaystyle\int\limits_{0^{+}}}
\rho^{-2}(u)du=+\infty$, there exist a decreasing sequence $(a_{n})$ of
positive real numbers such that $1>a_{1}$

satisfying

\begin{center}
$%
{\displaystyle\int\limits_{a_{1}}^{1}}
\rho^{-2}(u)du=1$, $%
{\displaystyle\int\limits_{a_{2}}^{a_{1}}}
\rho^{-2}(u)du=2,....,%
{\displaystyle\int\limits_{a_{n}}^{a_{n-1}}}
\rho^{-2}(u)du=n,.....$
\end{center}

Clearly $\left(  a_{n}\right)  $ converges to $0$ as $n$ tends to $+\infty$.

The properties of $\rho$ allow us to construct a sequence of functions
$\psi_{n}(u),$ $n=1,2,...$ , such that

i) $\psi_{n}(u)$ is a continuous function such that its support is contained
in $\left(  a_{n},a_{n-1}\right)  $

ii) $0\leq\psi_{n}(u)\leq\dfrac{2}{n}\rho^{-2}(u)$ and $%
{\displaystyle\int\limits_{a_{n}}^{a_{n-1}}}
\psi_{n}(u)du=1$

Let $\varphi_{n}(x)=%
{\displaystyle\int\limits_{0}^{\left\vert x\right\vert }}
dy%
{\displaystyle\int\limits_{0}^{y}}
\psi_{n}(u)du,$ $x\in\mathbb{R}$

It is clear that $\varphi_{n}\in\mathcal{C}^{2}(\mathbb{R})$ such that
$\left\vert \varphi_{n}^{\prime}\right\vert \leq1$ and $\left(  \varphi
_{n}\right)  $ is an increasing sequence converging to $\left\vert
x\right\vert .$

Let $X_{t}^{1}$ and $X_{t}^{2}$ two solutions of corresponding to the same
Brownian motion and the same MVSDE

\begin{center}
$X_{t}^{1}-X_{t}^{2}=%
{\displaystyle\int\limits_{0}^{t}}
\left(  \sigma(s,X_{s}^{1}\right)  -\sigma(s,X_{s}^{2}))dW_{s}+%
{\displaystyle\int\limits_{0}^{t}}
\left(  b(s,X_{s}^{1},\mathbb{P}_{X_{s}^{1}}\right)  -b(s,X_{s}^{2}%
,\mathbb{P}_{X_{s}^{2}}))dW_{s}$
\end{center}

By using It\^{o}'s formula we obtain

\begin{center}
$%
\begin{array}
[c]{cl}%
\varphi_{n}(X_{t}^{1}-X_{t}^{2})= &
{\displaystyle\int\limits_{0}^{t}}
\varphi_{n}^{\prime}(X_{s}^{1}-X_{s}^{2})\left(  \sigma(s,X_{s}^{1}\right)
-\sigma(s,X_{s}^{2}))dW_{s}\\
& +%
{\displaystyle\int\limits_{0}^{t}}
\varphi_{n}^{\prime}(X_{s}^{1}-X_{s}^{2})\left(  b(s,X_{s}^{1},\mathbb{P}%
_{X_{s}^{1}})-b(s,X_{s}^{2},\mathbb{P}_{X_{s}^{2}})\right)  ds\\
& +\dfrac{1}{2}%
{\displaystyle\int\limits_{0}^{t}}
\varphi_{n}^{\prime\prime}(X_{s}^{1}-X_{s}^{2})\left(  \sigma(s,X_{s}%
^{1})-\sigma(s,X_{s}^{2})\right)  ^{2}ds
\end{array}
$
\end{center}

$\varphi_{n}^{\prime}$ and $\sigma$ being bounded, then the process under the
sign integral is sufficiently integrable. Then the first term is a true
martingale, so that its expectation is $0.$ Therefore

$%
\begin{array}
[c]{cl}%
E\left(  \varphi_{n}(X_{t}^{1}-X_{t}^{2})\right)  = & E\left[
{\displaystyle\int\limits_{0}^{t}}
\varphi_{n}^{\prime}(X_{s}^{1}-X_{s}^{2})\left(  b(s,X_{s}^{1},\mathbb{P}%
_{X_{s}^{1}})-b(s,X_{s}^{2},\mathbb{P}_{X_{s}^{2}})\right)  ds\right] \\
& +\dfrac{1}{2}E\left[
{\displaystyle\int\limits_{0}^{t}}
\varphi_{n}^{\prime\prime}(X_{s}^{1}-X_{s}^{2})\left(  \sigma(s,X_{s}%
^{1})-\sigma(s,X_{s}^{2})\right)  ^{2}ds\right] \\
& =I_{1}+I_{2}%
\end{array}
$

But we know that $W_{1}(\mathbb{P}_{X_{s}^{1}},\mathbb{P}_{X_{s}^{2}%
})=E\left(  \left\vert X_{s}^{1}\mathbb{-}X_{s}^{2}\right\vert \right)  $

Then

\begin{center}
$\left\vert I_{1}\right\vert \leq E%
{\displaystyle\int\limits_{0}^{t}}
\kappa(\left\vert X_{s}^{1}\mathbb{-}X_{s}^{2}\right\vert )ds+%
{\displaystyle\int\limits_{0}^{t}}
CE\left(  \left\vert X_{s}^{1}\mathbb{-}X_{s}^{2}\right\vert \right)  ds$
\end{center}

Then by Growall lemma, there exist a constant $M$ such that $\left\vert
I_{1}\right\vert \leq M.E%
{\displaystyle\int\limits_{0}^{t}}
\kappa(\left\vert X_{s}^{1}\mathbb{-}X_{s}^{2}\right\vert )ds$

On the other hand

\begin{center}
$%
\begin{array}
[c]{cc}%
\left\vert I_{2}\right\vert = & \dfrac{1}{2}E\left[
{\displaystyle\int\limits_{0}^{t}}
\varphi_{n}^{\prime\prime}(X_{s}^{1}-X_{s}^{2})\left(  \sigma(s,X_{s}%
^{1})-\sigma(s,X_{s}^{2})\right)  ^{2}ds\right] \\
& \leq\dfrac{1}{2}E\left[
{\displaystyle\int\limits_{0}^{t}}
\dfrac{2}{n}\rho^{-2}(X_{s}^{1}-X_{s}^{2})\rho^{2}(X_{s}^{1}-X_{s}%
^{2})ds\right]  =\dfrac{t}{n}%
\end{array}
$
\end{center}

Then $\left\vert I_{2}\right\vert $ tends to $0$ as $n$ tends to $+\infty.$

Letting $n$ tending to $+\infty$ it holds that: $E\left(  \left\vert X_{t}%
^{1}-X_{t}^{2}\right\vert \right)  \leq M.E%
{\displaystyle\int\limits_{0}^{t}}
\kappa(\left\vert X_{s}^{1}\mathbb{-}X_{s}^{2}\right\vert )ds$. Since $%
{\displaystyle\int\limits_{0^{+}}}
\kappa^{-1}(u)du=+\infty$ $\ $we conclude that $E\left(  \left\vert X_{t}%
^{1}-X_{t}^{2}\right\vert \right)  =0.$\eop

\textbf{Remark. }\textit{The continuity and boundness of the coefficients
imply the existence of a weak solution (see \cite{JMW} Proposition 1.10 ).
Then by the well known Yamada - Watanabe theorem applied to equation
(\ref{MVSDE2}) (see \cite{Kur} example 2.14, page 10), the pathwise uniqueness
proved in the last theorem implies the existence and uniqueness of a strong
solution.}

\section{Convergence of the Picard successive \textbf{approximation }}

Assume that $b(t,x,\mu)$ and\textbf{ }$\sigma(t,x,\mu)$ satisfy assumptions
$\mathbf{(H_{1})}$, $\mathbf{(H_{2}).}$ We will prove the convergence of the
Picard iteration scheme. This scheme is useful for numerical computations of
the unique solution of (\ref{MVSDE}). Let $(X_{t}^{0})=x$ for all $t\in\left[
0,T\right]  $ and define $\left(  X_{t}^{n+1}\right)  $ as the solution of the
following SDE

\begin{center}
$\left\{
\begin{array}
[c]{c}%
dX_{t}^{n+1}=b(t,X_{t}^{n},\mathbb{P}_{X_{t}^{n}})dt+\sigma(t,X_{t}%
^{n},\mathbb{P}_{X_{t}^{n}})dB_{t}\\
X_{0}^{n+1}=x
\end{array}
\right.  $
\end{center}

\begin{theorem}
Under assumptions $\mathbf{(H_{1})}$, $\mathbf{(H_{2})}$, the sequence
$\left(  X^{n}\right)  $ converges to the unique solution of \ (\ref{MVSDE})
\[
E[\sup_{t\leq T}|X_{t}^{n}-X_{t}|^{2}]\rightarrow0
\]

\end{theorem}

\bop Let $n\geq0,$ by applying usual arguments such as Schwartz inequality and
Burkholder-Davis Gundy inequality for the martingale part, we get%

\begin{align*}
|X_{t}^{n+1}-X_{t}^{n}|^{2}  &  \leq2(\int_{0}^{t}|b(s,X_{s}^{n},P_{X_{s}^{n}%
})-b(s,X_{s}^{n-1},P_{X_{s}^{n-1}})|ds)^{2}\\
&  +2(\int_{0}^{t}|\sigma(s,X_{s}^{n},P_{X_{s}^{n}})-\sigma(s,X_{s}%
^{n-1},P_{X_{s}^{n-1}})|dB_{s})^{2}\\
E[\sup_{t\leq T}|X_{t}^{n+1}-X_{s}^{n}|^{2}]  &  \leq2TE[\int_{0}%
^{T}|b(s,X_{s}^{n},\mathbb{P}_{X_{s}^{n}})-b(s,X_{s}^{n-1},\mathbb{P}%
_{X_{s}^{n-1}})|^{2}ds]\\
&  +2C_{2}E[\int_{0}^{T}|\sigma(s,X_{s}^{n},\mathbb{P}_{X_{s}^{n}}%
)-\sigma(s,X_{s}^{n-1},\mathbb{P}_{X_{s}^{n-1}})|^{2}ds]
\end{align*}

the coefficients $b$ and $\sigma$ being Lipschitz continuous in $\left(
x,\mu\right)  $ we get%

\begin{align*}
E[\sup_{t\leq T}|X_{t}^{n+1}-X_{t}^{n}|^{2}]  &  \leq2(T+C_{2})L^{2}\int
_{0}^{T}E[|X_{s}^{n}-X_{s}^{n-1}|^{2}]+W_{2}(\mathbb{P}_{X_{s}^{n}}%
,\mathbb{P}_{X_{s}^{n-1}})ds\\
&  \leq4(T+C_{2})L^{2}\int_{0}^{T}E[|X_{s}^{n}-X_{s}^{n-1}|^{2}]ds\\
&  \leq4(T+C_{2})L^{2}\int_{0}^{T}E[\sup_{t\leq T}|X_{s}^{n}-X_{s}^{n-1}%
|^{2}]ds
\end{align*}

Then for all $n\geq1$, and $t\leq T$%

\begin{align*}
E[\sup_{t\leq T}|X_{t}^{1}-X_{s}^{0}|^{2}]  &  \leq2T\int_{0}^{T}%
b|(s,x,\mu)|^{2}ds+C_{2}\int_{0}^{T}\sigma|(s,x,\mu)|^{2}ds\\
&  \leq2(C_{2}+T)M(1+E(|x|^{2}))T\\
&  \leq A_{1}T
\end{align*}

where the constant $A_{1}$ only depends on $C_{2},M,T$ and $E[|x|^{2}]$. So by
induction on $n$ we obtain%

\[
E[\sup_{t\leq T}|X_{t}^{n+1}-X_{t}^{n}|^{2}]\leq\frac{A_{2}^{n+1}T^{n+1}%
}{(n+1)!}%
\]

This implies in particular that $\left(  X_{t}^{n}\right)  $ is a Cauchy
sequence in $L^{2}(\Omega,\mathcal{C}(\left[  0,T\right]  ,\mathbb{R}^{d}))$
which is complete. Therefore $\left(  X_{t}^{n}\right)  $ converges to a limit
$\left(  X_{t}\right)  $ which is the unique solution of (\ref{MVSDE}) \eop

\section{\textbf{Stability with respect to initial condition}}

In this section, we will study the stability of MFSDEs with respect to small
perturbations of the initial condition.

We denote by $\left(  X_{t}^{x}\right)  $ the unique solution of (\ref{MVSDE})
such that $X_{0}^{x}=x$%

\[
\left\{
\begin{array}
[c]{c}%
dX_{t}^{x}=b(t,X_{t}^{x},\mathbb{P}_{X_{t}^{x}})dt+\sigma(t,X_{t}%
^{x},\mathbb{P}_{X_{t}^{x}})dB_{t}\\
X_{0}^{x}=x
\end{array}
\right.
\]

\begin{theorem}
Assume that $b(t,x,\mu)$ and $\sigma(t,x,\mu)$ satisfy $\mathbf{(H_{1})}$,
$\mathbf{(H_{2}),}$ then the mapping
\end{theorem}

\begin{center}
$\Phi:\mathbb{R}^{d}\longrightarrow L^{2}(\Omega,\mathcal{C}(\left[
0,T\right]  ,\mathbb{R}^{d}))$
\end{center}

defined by $\left(  \Phi(x)_{t}\right)  =\left(  X_{t}^{x}\right)  $ is continuous.

\bop

Let $\left(  x_{n}\right)  $ be a sequence in $\mathbb{R}^{d}$ converging to
$x.$ Let us prove that $\underset{n\longrightarrow+\infty}{\lim}E\left[
\sup\limits_{t\leq T}|X_{t}^{n}-X_{t}|^{2}\right]  =0,$ where $X_{t}^{n}%
=X_{t}^{x_{n}}.$ We have%

\begin{align*}
|X_{t}^{n}-X_{t}|^{2}  &  =|x_{n}-x+\int_{0}^{t}(b(s,X_{s}^{n},\mathbb{P}%
_{X_{s}^{n}})-b(s,X_{s}^{n},\mathbb{P}_{X_{s}^{n}}))ds\\
&  +\int_{0}^{t}(\sigma(s,X_{s}^{n},\mathbb{P}_{X_{s}^{n}})-\sigma(s,X_{s}%
^{n},\mathbb{P}_{X_{s}^{n}}))dB_{s}|^{2}%
\end{align*}%
\begin{align*}
&  \leq3|x_{n}-x|^{2}+3(\int_{0}^{t}|b(s,X_{s}^{n},\mathbb{P}_{X_{s}^{n}%
})-b(s,X_{s}^{n},\mathbb{P}_{X_{s}^{n}})|ds)^{2}\\
&  +3(\int_{0}^{t}|\sigma(s,X_{s}^{n},\mathbb{P}_{X_{s}^{n}})-\sigma
(s,X_{s}^{n},\mathbb{P}_{X_{s}^{n}})|dB_{s})^{2}%
\end{align*}

\begin{align*}
E\left[  \sup_{t\leq T}|X_{t}^{n}-X_{t}|^{2}\right]   &  \leq3|x_{n}%
-x|^{2}+3E[\sup_{s\leq t}\int_{0}^{t}|b(s,X_{s}^{n},\mathbb{P}_{X_{s}^{n}%
})-b(s,X_{s}^{n},\mathbb{P}_{X_{s}^{n}})|ds]^{2}\\
&  +3E[\sup_{s\leq t}\int_{0}^{t}|\sigma(s,X_{s}^{n},\mathbb{P}_{X_{s}^{n}%
})-\sigma(s,X_{s}^{n},\mathbb{P}_{X_{s}^{n}})|dB_{s}]^{2}%
\end{align*}

we apply Schwartz and Burkholder Davis Gundy inequalities to obtain%

\begin{align*}
E\left[  \sup_{t\leq T}|X_{t}^{n}-X_{t}|^{2}\right]   &  \leq3|x_{n}%
-x|^{2}+3TE[\int_{0}^{t}|b(s,X_{s}^{n},\mathbb{P}_{X_{s}^{n}})-b(s,X_{s}%
^{n},\mathbb{P}_{X_{s}^{n}})|^{2}ds]\\
&  +3C_{2}E[\int_{0}^{t}|\sigma(s,X_{s}^{n},\mathbb{P}_{X_{s}^{n}}%
)-\sigma(s,X_{s}^{n},\mathbb{P}_{X_{s}^{n}})|^{2}ds]
\end{align*}

The Lipschitz condition implies that$\ $%
\[
E\left[  \sup_{t\leq T}|X_{t}^{n}-X_{t}|^{2}\right]  \leq3|x_{n}%
-x|^{2}+3(T+C_{2})L^{2}[\int_{0}^{t}E|X_{s}^{n}-X_{s}|^{2}+W_{2}%
(\mathbb{P}_{X_{s}^{n}},\mathbb{P}_{X_{s}})]ds
\]

Since%
\[
W_{2}^{2}(\mathbb{P}_{X_{t}^{n}},\mathbb{P}_{X_{t}})\leq E[|X_{s}^{n}%
-X_{s}|^{2}],
\]

then%

\begin{align*}
E\left[  \sup_{t\leq T}|X_{t}^{n}-X_{t}|^{2}\right]   &  \leq3|x_{n}%
-x|^{2}+6(T+c_{2})L^{2}\int_{0}^{t}E|X_{s}^{n}-X_{s}|^{2}ds\\
&  \leq3|x_{n}-x|^{2}+6(T+c_{2})L^{2}\int_{0}^{t}E[\sup_{t\leq T}|X_{s}%
^{n}-X_{s}|^{2}]ds.
\end{align*}

Finally we apply Gronwall lemma to conclude that%

\[
E\left[  \sup_{t\leq T}|X_{t}^{n}-X_{t}|^{2}\right]  \leq3|x_{n}-x|^{2}%
\exp[6(T+c_{2})L^{2}t]
\]

Therefore $\lim_{n\rightarrow\infty}x_{n}=x$ implies that $\lim_{n\rightarrow
\infty}E\left[  \sup_{t\leq T}|X_{t}^{n}-X_{t}|^{2}\right]  =0.$\eop

\section{\textbf{Stability with respect to the coefficients}}

In this section, we will establish the stability of the MVSDE with respect to
small perturbation of the cofficients $b$ and $\sigma.$ Let us consider
sequences of functions $\left(  b_{n}\right)  $ and $\left(  \sigma
_{n}\right)  $ and consider the corresponding MFSDE:%
\begin{align}
dX_{t}^{n}  &  =b_{n}(t,X_{t}^{n},\mathbb{P}_{X_{t}^{n}})dt+\sigma_{n}%
(t,X_{t}^{n},\mathbb{P}_{X_{t}^{n}})dB_{t}\label{MVSDEn}\\
X_{0}^{n}  &  =x\nonumber
\end{align}

The following theorem gives us the continuous dependence of the solution with
respect to the coefficients.

\begin{theorem}
Assume that the functions $b(t,x,\mu),$ $b_{n}(t,x,\mu),$ $\sigma(t,x,\mu)$
and $\sigma_{n}(t,x,\mu)$ satisfy $\mathbf{(H_{1})}$, $\mathbf{(H_{2}).}$
Further suppose that for each $T>0,$ and each compact set $K$ there existe
$C>0$ such that{}
\end{theorem}

$i)\sup_{t\leq T}(|b_{n}(t,x,\mu)|+|\sigma_{n}(t,x,\mu)|)\leq C(1+|x|),$

$ii)\lim_{n\rightarrow\infty}\sup_{t\leq T}\sup_{x\in K}\sup_{\mu
\in\mathcal{P}_{2}(\mathbb{R}^{d})}||b_{n}(t,x,\mu)-b(t,x,\mu)||+||\sigma
_{n}(t,x,\mu)-\sigma(t,x,\mu)||=0$

then%

\[
\lim_{n\rightarrow\infty}E\left[  \sup_{t\leq T}|X_{t}^{n}-X_{t}|^{2}\right]
=0
\]

where $\left(  X_{t}^{n}\right)  $ and $\left(  X_{t}\right)  $ are
respectively solutions of (\ref{MVSDEn}) and (\ref{MVSDE}).

\bop

For each $n\in\mathbb{N}$, let $\left(  X_{t}^{n}\right)  $ be a solution of
(\ref{MVSDEn}), then by using%

\begin{align*}
|X_{t}^{n}-X_{t}|^{2}  &  \leq3(\int_{0}^{t}|b_{n}(s,X_{s}^{n},\mathbb{P}%
_{X_{s}^{n}})-b_{n}(s,X_{s},\mathbb{P}_{X_{s}})|ds)^{2}\\
&  +3(\int_{0}^{t}|b_{n}(s,X_{s},\mathbb{P}_{X_{s}})-b(s,X_{s},\mathbb{P}%
_{X_{s}})|ds)^{2}\\
&  +3\left\vert \int_{0}^{t}\left(  \sigma_{n}(s,X_{s}^{n},\mathbb{P}%
_{X_{s}^{n}})-\sigma_{n}(s,X_{s},\mathbb{P}_{X_{s}})\right)  dB_{s}\right\vert
^{2}\\
&  +3\left\vert \int_{0}^{t}\left(  \sigma_{n}(s,X_{s},\mathbb{P}_{X_{s}%
})+\sigma(s,X_{s},\mathbb{P}_{X_{s}})\right)  dB_{s}\right\vert ^{2}%
\end{align*}

By using the Lipschitz continuity and Burkholder Davis Gundy inequality, it
holds that%

\begin{align*}
E\left[  \sup_{t\leq T}|X_{t}^{n}-X_{t}|^{2}\right]   &  \leq3(T+C_{2}%
)L^{2}\int_{0}^{t}E[|X_{s}^{n}-X_{s}|^{2}]+W_{2}(\mathbb{P}_{X_{s}^{n}%
},\mathbb{P}_{X_{s}})^{2}]ds\\
&  +3(T+C_{2})E[\int_{0}^{t}|b_{n}(s,X_{s},\mathbb{P}_{X_{s}})-b(s,X_{s}%
,\mathbb{P}_{X_{s}})|^{2}ds]\\
&  +3(T+C_{2})E[\int_{0}^{t}|\sigma_{n}(s,X_{s},\mathbb{P}_{X_{s}}%
)-\sigma(s,X_{s},\mathbb{P}_{X_{s}})|^{2}ds]\\
&  \leq6(T+C_{2})L^{2}\int_{0}^{T}E[|X_{s}^{n}-X_{s}|^{2}]ds+K_{n}\\
&  \leq6(T+C_{2})L^{2}\int_{0}^{T}E\left[  \sup_{s\leq t}|X_{s}^{n}-X_{s}%
|^{2}\right]  dt+K_{n}%
\end{align*}

such that
\[
K_{n}=3(T+C_{2})E[\int_{0}^{T}\left(  |b_{n}(s,X_{s},\mathbb{P}_{X_{s}%
})-b(s,X_{s},\mathbb{P}_{X_{s}})|^{2}+|\sigma_{n}(s,X_{s},\mathbb{P}_{X_{s}%
})+\sigma(s,X_{s},\mathbb{P}_{X_{s}})|^{2}\right)  ds]
\]

An application of Gronwall lemma allows us to get%

\[
E\left[  \sup_{t\leq T}|X_{t}^{n}-X_{t}|^{2}\right]  \leq K_{n}\exp
6(T+C_{2})L^{2}.T
\]

By using assumptions i) and ii) it is easy to see that $K_{n}$
$\longrightarrow0.$as $n\longrightarrow+\infty,$which achieves the
proof.\eop

\section{Stability\textbf{ with respect to the driving processes}}

In this section, we consider McKean-Vlasov SDE driven by continuous semi-martingales.

Let \ $b:[0,T]\times\mathbb{R}^{d}\times\mathcal{P}_{2}(\mathbb{R}%
^{d})\rightarrow\mathbb{R}^{d}$ and $\sigma:[0,T]\times\mathbb{R}^{d}%
\times\mathcal{P}_{2}(\mathbb{R}^{d})\rightarrow\mathbb{R}^{d\times d}$ be
bounded continuous functions.

We consider MVSDEs driven by continuous semi-martingales of the following
form
\begin{equation}
\left\{
\begin{array}
[c]{c}%
dX_{t}=b(t,X_{t},\mathbb{P}_{X_{t}})dA_{t}+\sigma(t,X_{t},\mathbb{P}_{X_{t}%
})dM_{t}\\
X_{0}=x
\end{array}
\right.  \label{MVSDE3}%
\end{equation}

where $A_{t}$ is an adapted continuous process of bounded variation and
$M_{t}$ is a continuous local martingale.

Let us consider the following sequence of MVSDEs%

\begin{equation}
\left\{
\begin{array}
[c]{c}%
dX_{t}^{n}=b(t,X_{t}^{n},\mathbb{P}_{X_{t}^{n}})dA_{t}^{n}+\sigma(t,X_{t}%
^{n},\mathbb{P}_{X_{t}^{n}})dM_{t}^{n}\\
X_{0}^{n}=x
\end{array}
\right.  \label{MVSDE3n}%
\end{equation}

where $(A^{n})$ is a sequence of $%
\mathcal{F}%
_{t}$-adapted continuous process of bounded variaton and $M^{n}$ is continuous
$(%
\mathcal{F}%
_{t},\mathbb{P)}$-local martingales.

Let us assume that $(A,A^{n},M,M^{n})$ satisfy:

(\textbf{H}$_{\mathbf{7}}$)

\qquad1) The family $(A,A^{n},M,M^{n})$ is bounded in $%
\mathbb{C}
([0,1])^{4}.$

\qquad2) $\left(  M^{n}-M\right)  $ converges to $0$ in probability in $%
\mathbb{C}
([0,1])$ as n tends to $+\infty.$

\qquad3)The total variation $(A^{n}-A)$ converges to $0$ in probability as n
tends to $+\infty.$

\begin{theorem}
Let $b(t,x,\mu)$ and $\sigma(t,x,\mu)$ satisfy $\mathbf{(H_{1})}$,
$\mathbf{(H_{2})}$. Further assume that $(A,A^{n},M,M^{n})$ satisfy
($\mathbf{H}_{\mathbf{7}}$). Then for each $\varepsilon>0$%
\[
\lim_{n\rightarrow\infty}E[\sup_{t\leq T}|X_{t}^{n}-X_{t}|^{2}]=0
\]

where $\left(  X_{t}^{n}\right)  $ and $\left(  X_{t}\right)  $ are
respectively solutions of (\ref{MVSDE3n}) and (\ref{MVSDE3}).
\end{theorem}

\bop

Let $n\in\mathbb{N}$, then by using similar arguments as in the preceding
theorems, we have%

\begin{align*}
\mathbb{E}[\sup_{t\leq T}|X_{t}^{n}-X_{t}|^{2}]  &  \leq3(E[\sup_{t\leq T}%
\int_{0}^{t}|b(s,X_{s}^{n},\mathbb{P}_{X_{s}^{n}})-b(s,X_{s}^{n}%
,\mathbb{P}_{X_{s}^{n}})|]dA_{s}^{n})^{2}\\
&  +3(E\left[  \sup_{t\leq T}\left\vert \int_{0}^{t}\left(  \sigma(s,X_{s}%
^{n},\mathbb{P}_{X_{s}^{n}})-\sigma(s,X_{s}^{n},\mathbb{P}_{X_{s}^{n}%
})\right)  dM_{s}^{n}\right\vert ^{2}\right] \\
&  +3E[(\sup_{t\leq T}\int_{0}^{t}\left\vert b(t,X_{t},\mathbb{P}_{X_{t}%
})\right\vert d\left\vert A_{s}^{n}-A_{s}\right\vert )^{2}+\sup_{t\leq
T}\left\vert \int_{0}^{t}\sigma(t,X_{t},\mathbb{P}_{X_{s}})d(M_{s}^{n}%
-M_{s})\right\vert ^{2}]
\end{align*}

Let%
\[
K_{n}=3E[(\sup_{t\leq T}\int_{0}^{t}\left\vert b(t,X_{t},\mathbb{P}_{X_{t}%
})\right\vert d\left\vert A_{s}^{n}-A_{s}\right\vert )^{2}+\sup_{t\leq
T}\left\vert \int_{0}^{t}\sigma(t,X_{t},\mathbb{P}_{X_{s}})d(M_{s}^{n}%
-M_{s})\right\vert ^{2}]
\]

By using Schwartz and Burkholder Davis Gundy inequalities along with the
Lipschitz condition, we obtain%

\begin{align*}
\mathbb{E}[\sup_{t\leq T}|X_{t}^{n}-X_{t}|^{2}]  &  \leq C(T)[\int_{0}%
^{T}(E\left(  \sup_{s\leq t}|X_{s}^{n}-X_{s}|^{2}\right)  +\mathbb{W}%
_{2}(\mathbb{P}_{X_{s}^{n}},\mathbb{P}_{X_{s}})^{2})(dA_{s}^{n}+d<M^{n}%
,M^{n}>_{s})]+K_{n}\\
&  \leq2C(T)\int_{0}^{T}E[\sup_{s\leq t}|X_{s}^{n}-X_{s}|^{2}](dA_{s}%
^{n}+d<M^{n},M^{n}>_{s})+K_{n}%
\end{align*}

where $C(T)$ is a positive constant which may change from line to line$.$

Since $(A_{s}^{n}+d<M^{n},M^{n}>_{s})$ is an increasing process, then
according to the Stochastic Gronwall lemma \cite{Met} Lemma 29.1, page 202, we have%

\[
\mathbb{E}[\sup_{t\leq T}|X_{t}^{n}-X_{t}|^{2}]\leq2K_{n}CE(A_{T}^{n}%
+<M^{n},M^{n}>_{T}))<+\infty,
\]

where $C$ is a constant.

By using assumption $\left(  \mathbf{H}_{\mathbf{7}}\right)  $ it is easy to
that
\[
\lim_{n\rightarrow\infty}K_{n}=0
\]

Therefore%

\[
\lim_{n\rightarrow\infty}\mathbb{E}[\sup_{t\leq T}|X_{t}^{n}-X_{t}|^{2}]=0
\]
\eop

\section{Approximation of relaxed control problems}

It is well known that in the deterministic as well as in stochastic control
problems , an optimal control does not necessarily exist\ in the space of
strict controls, in the absence of convexity conditions. The classical method
is then to introduce measure valued controls which describe the introduction
of a stochastic parameter see \cite{EKNJ} and the references therein. These
measure valued controls called relaxed controls generalize the strict controls
in the sense that the set of strict controls may be identified as a dense
subset of the set of the relaxed controls. The relaxed control problem is a
true extension of the strict control problem if they have the same value
function.\ That is the infimum among strict controls is equal to the infimum
among relaxed controls. This last property is based on the continuity of the
dynamics and the cost functional with respect to the control variable. We show
that under Lipschitz condition and continuity with respect to the control
variable of the coefficients that the strict and relaxed control problems have
the same value function. Our result extends those in \cite{BMM1, BMM2}, to
general MFSDEs of the type \ref{CMFSDE}.

Let $\mathbb{A}$ be some compact metric space called the action space. A
strict control $\left(  u_{t}\right)  $ is a measurable, $\mathcal{F}_{t}-$
adapted process with values in the action space $\mathbb{A}$. We denote
$\mathcal{U}_{ad}$ the space of strict controls.

\noindent The state process corresponding to a strict control is the unique
solution, of the following MFSDE%

\begin{equation}
\left\{
\begin{array}
[c]{l}%
dX_{t}=b(t,X_{t},\mathbb{P}_{X_{t}},u_{t})ds+\sigma(t,X_{t},\mathbb{P}_{X_{t}%
},u_{t})dB_{s}\\
X_{0}=x
\end{array}
\right.  \label{CMFSDE}%
\end{equation}

\noindent and the corresponding cost functional is given by

\begin{center}
$J(u)=E\left[  \int_{0}^{T}h(t,X_{t},\mathbb{P}_{X_{t}},u_{t}dt+g(X_{T}%
,\mathbb{P}_{X_{T}})\right]  .$
\end{center}

\noindent The problem is to minimize $J(u)$ over the space $\mathcal{U}_{ad}$
of strict controls and to find $u^{\ast}\in$ $\mathcal{U}_{ad}$ such that
$J(u^{\ast})=\inf\left\{  J(u),u\in\mathcal{U}_{ad}\right\}  .$

\noindent Let us consider the following assumptions in this section.

$\mathbf{(H}_{\mathbf{4}}\mathbf{)}$ $b:\left[  0,T\right]  \times
\mathbb{R}^{d}\times\mathcal{P}_{2}(\mathbb{R}^{d})\times\mathbb{A}%
\longrightarrow\mathbb{R}^{d}$, $\sigma:\left[  0,T\right]  \times
\mathbb{R}^{d}\times\mathcal{P}_{2}(\mathbb{R}^{d})\times\mathbb{A}%
\longrightarrow\mathbb{R}^{d}\otimes\mathbb{R}^{d},$ are continuous bounded
functions $.$\ 

$\mathbf{(H}_{\mathbf{5}}\mathbf{)}$ $b(t,.,.,a)$ and $\sigma(t,.,.,a)$ are
Lipschitz continuyous uniformly in $(t,a)\in\left[  0,T\right]  \times
\mathbb{A}.$

$\mathbf{(H}_{\mathbf{6}}\mathbf{)}$ $h:\left[  0,T\right]  \times
\mathbb{R}\times\mathbb{R}\times\mathbb{A}\longrightarrow\mathbb{R}$ and
$g:\mathbb{R}\times\mathbb{R}\longrightarrow\mathbb{R}$, are bounded
continuous functions, such that $h(t,.,.,a)$ is Lipschitz in $(x,\mu).$

\noindent It is clear that under assumptions $\mathbf{(H}_{\mathbf{4}%
}\mathbf{)}$ and $\mathbf{(H}_{\mathbf{5}}\mathbf{)}$ and according to Theorem
3.1 , for each $u\in\mathcal{U}_{ad},$ the MFSDE (\ref{CMFSDE}) has a unique
strong solution, such that for every $p>0$, $E(\left\vert X_{t}\right\vert
^{p})<+\infty.$ Moreover for each $u\in$ $\mathcal{U}_{ad}$ $\left\vert
J(u)\right\vert <+\infty$.

Let $\mathbb{V}\ $be the set of product measures $\mu$ on $[0,T]\times
\mathbb{A}$ whose projection on $\left[  0,T\right]  $ coincides with the
Lebesgue measure $dt$. $\mathbb{V}$ as a closed subspace of the space of
positive Radon measures $\mathbb{M}_{+}([0,T]\times\mathbb{A})$ is compact for
the topology of weak convergence.

\begin{definition}
A relaxed control on the filtered probability space $\left(  \Omega
,\mathcal{F},\mathcal{F}_{t},P\right)  $ is a random variable $\mu=dt.\mu
_{t}(da)$ with values in $\mathbb{V}$, such that $\mu_{t}(da)$ is
progressively measurable with respect to $(\mathcal{F}_{t})$ and such that for
each $t$, $1_{(0,t]}.\mu$ is $\mathcal{F}_{t}-$measurable.
\end{definition}

\begin{remark}
The set $\mathcal{U}_{ad}$ of strict controls is embedded into the set of
relaxed controls by identifying $u_{t}$ with $dt\delta_{u_{t}}(da).$
\end{remark}

It was proved in \cite{EM} for classical control problems and in \cite{BMM2}
that the relaxed state process corresponding to a relaxed control must satisfy
a MFSDE driven by a martingale measure instead of a Brownian motion. That is
the relaxed state process satisfies%

\begin{equation}
\left\{
\begin{array}
[c]{c}%
dX_{t}=%
{\textstyle\int_{\mathbb{A}}}
b(t,X_{t},\mathbb{P}_{X_{t}},a)\,\mu_{t}(da)dt+%
{\textstyle\int_{\mathbb{A}}}
\sigma(t,X_{t},\mathbb{P}_{X_{t}},a)\,M(da,dt)\\
X_{0}=x,
\end{array}
\right.  \label{REMFSDE}%
\end{equation}

where $M\ $is an orthogonal continuous martingale measure, with
intensity\textit{ }$dt\mu_{t}(da).$Using the same tools as in Theorem 3.1, it
is not difficult to prove that (\ref{REMFSDE}) admits a unique strong
solution. The following Lemma, known in the control literature as Chattering
Lemma states that the set of strict controls is a dense subset in the set of
relaxed controls.

\begin{lemma}
i)\textbf{ }\textit{Let }$(\mu_{t})$\textit{\ be a relaxed control}%
$.$\textit{\ Then there exists a sequence of adapted processes }$(u_{t}^{n}%
)$\textit{\ with values in }$\mathbb{A}$\textit{, such that the sequence of
random measures }$\left(  \delta_{u_{t}^{n}}(da)\,dt\right)  $%
\textit{\ converges in }$\mathbb{V}$\textit{ to }$\mu_{t}(da)\,dt,$\textit{
}$P-a.s.$

ii) For any $g$ continuous in $\left[  0,T\right]  \times\mathbb{M}%
_{1}(\mathbb{A})$ such that $g(t,.)$ is linear, we have

$\underset{n\rightarrow+\infty}{\lim}\int_{0}^{t}g(s,\delta_{u_{s}^{n}%
})ds=\int_{0}^{t}g(s,\mu_{s})ds$ uniformly in $t\in\left[  0,T\right]  ,$
$P-a.s.$
\end{lemma}

\bop See \cite{EKNJ} \eop

\textit{Let }$X_{t}^{n}$ be the solution \ of the state equation\textit{(}%
\ref{CMFSDE}) corresponding to $u^{n},$ where $u^{n}$ is a strict control
defined as in the last Lemma. If we denote $M^{n}(t,F)=\int\nolimits_{0}%
^{t}\int\nolimits_{F}\delta_{u_{s}^{n}}(da)dW_{s},$ then $M^{n}(t,F)$ is an
orthogonal martingale measure and $X_{t}^{n}$ may be written in a relaxed form
as follows

\begin{center}
$\left\{
\begin{array}
[c]{l}%
dX_{t}^{n}=%
{\displaystyle\int\limits_{\mathbb{A}}}
b(t,X_{t}^{n},\mathbb{P}_{X_{t}^{n}},a)\delta_{u_{t}^{n}}(da)dt+%
{\displaystyle\int\limits_{\mathbb{A}}}
\sigma(t,X_{t},\mathbb{P}_{X_{t}^{n}},a)M^{n}(dt,da)\\
X_{0}=x
\end{array}
\right.  $
\end{center}

Therefore $X_{t}^{n}$ may be viewed as the solution of \textit{(}%
\ref{REMFSDE}) corresponding to the relaxed control $\mu^{n}=dt\delta
_{u_{t}^{n}}(da).$

Since $\left(  \delta_{u_{t}^{n}}(da)\,dt\right)  $\textit{\ converges weakly
to }$\mu_{t}(da)\,dt,$\textit{ }$P-a.s.,$ then for every bounded predictable
process $\varphi:\Omega\times\left[  0,T\right]  \times\mathbb{A}%
\rightarrow\mathbb{R}$, such that $\varphi(\omega,t,.)$ is continuous$,$ we
have%
\begin{equation}
E\left[  \left(  \int\nolimits_{0}^{T}\int\nolimits_{\mathbb{A}}\varphi
(\omega,t,a)M^{n}(dt,da)-\int\nolimits_{0}^{t}\int\nolimits_{\mathbb{A}%
}\varphi(\omega,t,a)M(dt,da)\right)  ^{2}\right]  \rightarrow0\text{ }as\text{
}n\longrightarrow+\infty. \label{martin-measure}%
\end{equation}

see (\cite{BDM, Mel}).

The following proposition gives the continuity of the dynamics \textit{(}%
\ref{REMFSDE}) with respect to the control variable.

\begin{proposition}
i) \textit{If }$X_{t},$ $X_{t}^{n}$ \textit{denote the solutions of state
equation (}\ref{REMFSDE}) corresponding to $\mu$ and $\mu^{n},$ \textit{then
\ }$\overset{}{\text{For each}}t\leq T,$ $\underset{n\rightarrow+\infty}{\lim
}E(\left\vert X_{t}^{n}-X_{t}\right\vert ^{2})=0.$

ii) Let $J(u^{n})$ and $J(\mu)$ be the expected costs corresponding
respectively to $u^{n}$ and $\mu,$ then $\left(  J\left(  u^{n}\right)
\right)  $ converges to $J\left(  \mu\right)  .$
\end{proposition}

\bop1) Let $X_{t},$ $X_{t}^{n}$ the solutions of the MVSDE \textit{(}%
\ref{REMFSDE}) corresponding to $\mu$ and $u^{n}$. \ We have%

\[%
\begin{array}
[c]{cl}%
\left\vert X_{t}-X_{t}^{n}\right\vert  & \leq\left\vert \int\nolimits_{0}%
^{t}\int\nolimits_{\mathbb{A}}b\left(  s,X_{s},\mathbb{P}_{X_{t}},u\right)
\mu_{s}(da).ds-\int\nolimits_{0}^{t}\int\nolimits_{\mathbb{A}}b\left(
s,X_{s}^{n},\mathbb{P}_{X_{t}^{n}},u\right)  \delta_{u_{s}^{n}}%
(da)ds\right\vert \\
& +\left\vert \int\nolimits_{0}^{t}\int\nolimits_{\mathbb{A}}\sigma\left(
s,X_{s},\mathbb{P}_{X_{t}},a\right)  M(ds,da)-\int\nolimits_{0}^{t}%
\int\nolimits_{\mathbb{A}}\sigma\left(  s,X_{s}^{n},\mathbb{P}_{X_{t}^{n}%
},a\right)  M^{n}(ds,da)\right\vert \\
& \leq\left\vert \int\nolimits_{0}^{t}\int\nolimits_{\mathbb{A}}b\left(
s,X_{s},\mathbb{P}_{X_{t}},u\right)  \mu_{s}(da).ds-\int\nolimits_{0}^{t}%
\int\nolimits_{\mathbb{A}}b\left(  s,X_{s},\mathbb{P}_{X_{t}},a\right)
\delta_{u_{s}^{n}}(da)ds\right\vert \\
& +\left\vert \int\nolimits_{0}^{t}\int\nolimits_{\mathbb{A}}b\left(
s,X_{s},\mathbb{P}_{X_{t}},u\right)  \delta_{u_{s}^{n}}(da).ds-\int
\nolimits_{0}^{t}\int\nolimits_{\mathbb{A}}b\left(  s,X_{s}^{n},\mathbb{P}%
_{X_{t}},a\right)  \delta_{u_{s}^{n}}(da)ds\right\vert \\
& +\left\vert \int\nolimits_{0}^{s}\int\nolimits_{\mathbb{A}}\sigma\left(
v,X_{v},\mathbb{P}_{X_{v}},a\right)  M(dv,da)-\int\nolimits_{0}^{t}%
\int\nolimits_{\mathbb{A}}\sigma\left(  v,X_{v},\mathbb{P}_{X_{v}},a\right)
M^{n}(dv,da)\right\vert \\
& +\left\vert \int\nolimits_{0}^{s}\int\nolimits_{\mathbb{A}}\sigma\left(
v,X_{v},\mathbb{P}_{X_{v}},a\right)  M^{n}(dv,da)-\int\nolimits_{0}^{t}%
\int\nolimits_{\mathbb{A}}\sigma\left(  v,X_{v}^{n},\mathbb{P}_{X_{v}^{n}%
},a\right)  M^{n}(dv,da)\right\vert
\end{array}
\]

Then by using Burkholder-Davis-Gundy inequality for the martingale part and
the fact that all the functions in equation \textit{(}\ref{REMFSDE}) are
Lipschitz continuous, it holds that

\begin{center}%
\[
E\left(  \left\vert X_{t}-X_{t}^{n}\right\vert ^{2}\right)  \leq
C\int\nolimits_{0}^{T}E\left(  \left\vert X_{s}-X_{s}^{n}\right\vert
^{2}+\mathbb{W}_{2}(\mathbb{P}_{X_{s}^{n}},\mathbb{P}_{X_{s}})^{2}\right)
dt+K_{n},
\]

\end{center}

where $C$ is a nonnegative constant and

\begin{center}
$%
\begin{array}
[c]{c}%
K_{n}=E\left(  \left\vert \int\nolimits_{0}^{t}\int\nolimits_{\mathbb{A}%
}b\left(  s,X_{s},\mathbb{P}_{X_{t}},u\right)  \mu_{s}(da)ds-\int
\nolimits_{0}^{t}\int\nolimits_{\mathbb{A}}b\left(  s,X_{s},\mathbb{P}_{X_{t}%
},a\right)  \delta_{u_{s}^{n}}(da)ds\right\vert ^{2}\right) \\
+E\left(  \left\vert \int\nolimits_{0}^{t}\int\nolimits_{\mathbb{A}}%
\sigma\left(  s,X_{s},\mathbb{P}_{X_{t}},a\right)  M(ds,da)-\int
\nolimits_{0}^{t}\int\nolimits_{\mathbb{A}}\sigma\left(  s,X_{s}%
,\mathbb{P}_{X_{t}},a\right)  M^{n}(ds,da)\right\vert ^{2}\right) \\
=I_{n}+J_{n}%
\end{array}
$
\end{center}

Using the fact that

\begin{center}
$\mathbb{W}_{2}(\mathbb{P}_{X_{s}^{n}},\mathbb{P}_{X_{s}})^{2}\leq E\left(
\left\vert X_{s}-X_{s}^{n}\right\vert ^{2}\right)  $,
\end{center}

we get

\begin{center}%
\begin{equation}
E\left(  \left\vert X_{t}-X_{t}^{n}\right\vert ^{2}\right)  \leq
2C\int\nolimits_{0}^{T}E\left(  \left\vert X_{s}-X_{s}^{n}\right\vert
^{2}\right)  dt+K_{n}.
\end{equation}

\end{center}

Since the sequence $\left(  \delta_{u_{t}^{n}}(da)\,dt\right)  $%
\textit{\ converges weakly to }$\mu_{t}(da)\,dt,$\textit{ }$P-a.s.$ and $b$ is
bounded and continuous in the control variable, then by applying the Lebesgue
dominated convergence theorem we get $\underset{n\rightarrow+\infty}{\lim
}I_{n}=0$. On the other hand since $\sigma$ is bounded and continuous in $a$,
applying (\ref{martin-measure}) we get $\underset{n\rightarrow+\infty}{\lim
}J_{n}=0.$ We conclude by using Gronwall's Lemma.

ii) Let $u^{n}$ and $\mu$ as in i) then%

\[%
\begin{array}
[c]{cl}%
\left\vert J\left(  u^{n}\right)  -J\left(  \mu\right)  \right\vert  & \leq
E\left[  \int\limits_{0}^{T}\int_{\mathbb{A}}\left\vert h(t,X_{t}%
^{n},\mathbb{P}_{X_{t}^{n}},a)-\,h(t,X_{t},\mathbb{P}_{X_{t}},a)\right\vert
\delta_{u_{t}^{n}}(da)\,dt\right] \\
& +E\left[  \left\vert \int\limits_{0}^{T}\int_{\mathbb{A}}h(t,X_{t}%
,\mathbb{P}_{X_{t}},a)\delta_{u_{t}^{n}}(da)\,dt-\int\limits_{0}^{T}%
\int_{\mathbb{A}}h(t,X_{t},\mathbb{P}_{X_{t}},a)\mu_{t}(da)\,dt\right\vert
\right] \\
& +E\left[  \left\vert g(X_{T}^{n},\mathbb{P}_{X_{T}^{n}})-g(X_{T}%
,\mathbb{P}_{X_{T}})\right\vert \right]
\end{array}
\]

The first assertion implies that the sequence $\left(  X_{t}^{n}\right)  $
converges to $X_{t}$ in probability$,$ then by using the assumptions on the
coeffcients $h$ and $g$ and the dominated convergence theorem it is easy to
conclude .

\eop

\begin{remark}
According to the last Proposition, it is clear that the infimum among relaxed
controls is equal to the infimum among strict controls, which implies the
value functions for the relaxed and strict models are the same.
\end{remark}

\end{document}